\documentclass{amsart}

\usepackage{amssymb,amsmath,amsthm,latexsym,stmaryrd}
\usepackage[pagebackref=false,colorlinks=true]{hyperref}
\usepackage{color}

\newtheorem{theorem}{Theorem}[section]
\newtheorem{proposition}{Proposition}[section]

\newtheorem{definition}{Definition}[section]

\newcommand{\ie}{i.\hspace{.5pt}e.\ }
\newcommand{\f}{\phi}
\newcommand{\g}{\tilde{g}}
\newcommand{\n}{\nabla}
\newcommand{\M}{(\mathcal{M},\A\f,\A\xi,\A\eta,\A{}g)}

\newcommand{\R}{\mathbb R}
\newcommand{\X}{\mathfrak X}

\newcommand{\LL}{\mathcal{L}}
\newcommand{\ttt}{\tilde\tau}
\newcommand{\DD}{\partial}

\newcommand{\lm}{\lambda}
\newcommand{\al}{\alpha}

\newcommand{\A}{\allowbreak{}}
\newcommand{\D}{\mathrm{d}\hspace{-0.5pt}}
\newcommand{\thmref}[1]{Theorem~\ref{#1}}

\newcommand{\propref}[1]{Proposition~\ref{#1}}

\DeclareMathOperator{\grad}{grad} % the gradient
\DeclareMathOperator{\Hess}{Hess} % the gradient
\DeclareMathOperator{\tr}{tr} % the trace
\begin{document}

\vspace{2cm}

\title[Gradient almost para-Ricci-like solitons on \dots]
{Gradient almost para-Ricci-like solitons on para-Sasaki-like Riemannian $\Pi$-manifolds}

\author{Hristo Manev}
\address{Medical University of Plovdiv, Faculty of Pharmacy,
Department of Medical Physics and Biophysics,   15-A Vasil Aprilov
Blvd.,   Plovdiv 4002,   Bulgaria;}
\email{hristo.manev@mu-plovdiv.bg}


%\date{January 1, 2007 and, in revised form, June 22, 2007.}

%\dedicatory{{\rm Dedicated to my doctoral advisor Prof. Dimitar Mekerov on the occasion of his 70th birthday}}
%\newcommand{\FilSupport}{Research supported by the project NI13-FMI-002 of the Scientific
%Research Fund at the University of Plovdiv}


\begin{abstract}
Gradient almost para-Ricci-like soli\-tons on para-Sasaki-like Riemannian $\Pi$-manifolds are studied. It is proved that these objects have constant soliton coefficients. For the soliton under study is shown that the corresponding scalar curvatures of the considered both metrics are equal and constant and its Ricci tensor is a constant multiple of the vertical component. Explicit example of a 3-dimensional para-Sasaki-like Riemannian $\Pi$-manifold is provided in support of the proved assertions.
\thanks{The author was partially supported by Project MU21-FMI-008 of the Scientific Research Fund, University of Plovdiv, Bulgaria}
\end{abstract}
\maketitle

\section{Introduction}\label{sect-1}

The concept of Ricci solitons is introduced in \cite{Ham82} by R.S. Hamilton, \ie a special solution of the Ricci flow equation.
A pseudo-Riemannian manifold $(M,g)$ admits a \emph{Ricci soliton} when its Ricci tensor $\rho$ has the form \cite{Ham82}
\begin{equation*}\label{Rs}
\begin{array}{l}
\rho = -\frac12 \mathcal{L}_v g - \lm\, g,
\end{array}
\end{equation*}
where $\mathcal{L}$ denotes the Lie derivative, $v$ is a vector field and $\lm$ is a constant.
After its introduction a detailed study on Riemannian Ricci solitons is made in \cite{Cao}.
The start of the investigation of Ricci solitons in contact Riemannian geometry is given with \cite{Shar}.

In \cite{ChoKim}, a generalization of the Ricci soliton, called \emph{$\eta$-Ricci soliton}, is introduced in the geometry of almost contact metric manifolds by:
\begin{equation*}\label{eRs}
\begin{array}{l}
\rho = -\frac12 \mathcal{L}_v g  - \lm\, g  - \nu\, \eta\otimes\eta,
\end{array}
\end{equation*}
where $\nu$ is also a constant. This work is followed by different studies on almost contact metric manifolds with various types of additional structures (e.g. \cite{IngBag,NagPre}).

In next years this concept and its generalizations became objects of interest in different fields: in paracontact geometry \cite{Bla15,PraHad}; in pseudo-Riemannian geometry \cite{BagIng12,Bro-etal,MM-Sol1}; in mathematical physics \cite{ChaVal96,Fri85}.

We study the geometry of almost paracontact almost paracomplex Riemannian manifolds, called here briefly Riemannian $\Pi$-manifolds. Their induced almost product structure on the paracontact distribution is traceless and the restriction on the paracontact distribution of the almost paracontact structure is an almost paracomplex structure. The start of the investigation of the Riemannian $\Pi$-manifolds is given in \cite{ManSta}, where they are called almost paracontact Riemannian manifolds of type $(n,n)$. Due to the presence of two associated metrics $g$ and $\g$ on these manifolds, in \cite{HM17}, we started the investigation of the so-called \emph{para-Ricci-like solitons}, which is a natural generalization of the $\eta$-Ricci soliton.

In the present work we introduce and study \emph{gradient almost para-Ricci-like solitons} on para-Sasaki-like Riemannian $\Pi$-manifolds.

The paper is structured as follows.
In Section 1, we present some introductory words about the problem under study.
In the next Section 2, we present to the readers' attention the main results of the paper. In Section 3, we recall preliminary background facts about the considered geometry.
Section 4 is devoted to the proof of the main results. We prove that the corresponding scalar curvatures of the considered both metrics are equal and constant. Moreover, the Ricci tensor is a constant multiple of the vertical component. In the final Section 4, we support the assertions by an explicit 3-dimensional example.

\section{Main results}

\begin{theorem}\label{thm:grad}
Let $\M$ be a $(2n+1)$-dimensional para-Sasaki-like Riemannian $\Pi$-manifold.
If it admits a gradient almost para-Ricci-like soliton with functions $(\lm,\mu,\nu)$ and a potential function $f$, then $\M$ has constant scalar curvatures regarding $g$ and $\g$
\[
\tau=\ttt=-2n
\]
and its Ricci tensor has the form
\[
\rho=-2n\,\eta\otimes\eta.
\]
\end{theorem}

\section{Preliminaries}\label{sect-2}
\subsection{Para-Sasaki-like Riemannian $\Pi$-Manifolds}

Let $\M$ be a \emph{Riemannian $\Pi$-manifold}, namely, $\mathcal{M}$ is an odd-dimensional differentiable manifold, equipped with a Rie\-mannian metric $g$ and a Riemannian $\Pi$-structure $(\f,\xi,\eta)$ consisting of a (1,1)-tensor field $\f$, a Reeb vector field $\xi$ and its dual 1-form $\eta$. The following basic identities and their immediately derived properties are valid:
\begin{equation*}\label{strM}
\begin{array}{c}
\f\xi = 0,\qquad \f^2 = I - \eta \otimes \xi,\qquad
\eta\circ\f=0,\qquad \eta(\xi)=1,\\ \
\tr \f=0,\qquad g(\f x, \f y) = g(x,y) - \eta(x)\eta(y),
\end{array}
\end{equation*}
\begin{equation*}\label{strM2}
\begin{array}{ll}
g(\f x, y) = g(x,\f y),\qquad &g(x, \xi) = \eta(x),
\\
g(\xi, \xi) = 1,\qquad &\eta(\n_x \xi) = 0,
\end{array}
\end{equation*}
where $I$ and $\n$ denote the identity transformation on $T\mathcal{M}$ and the Levi--Civita connection of $g$, respectively (\cite{Sato76,ManTav57}).
Here and further, $x$, $y$, $z$, $w$ stand for arbitrary vector fields from $\X(\mathcal{M})$ or vectors in $T\mathcal{M}$ at a fixed point of $\mathcal{M}$.

Through $\g(x,y)=g(x,\f y)+\eta(x)\eta(y)$ is defined the associated metric $\g$ of $g$ on $\M$, which is an indefinite metric of signature $(n + 1, n)$ and compatible with the manifold in the same way as $g$.

In \cite{IvMaMa2}, a subclass called \emph{para-Sasaki-like Riemannian $\Pi$-manifolds} of the considered manifolds is introduced and studied. It is determined as follows
\begin{equation*}\label{defSl}
\begin{array}{l}
\left(\nabla_x\f\right)y=-g(x,y)\xi-\eta(y)x+2\eta(x)\eta(y)\xi\\
\phantom{\left(\nabla_x\f\right)y}=-g(\f x,\f y)\xi-\eta(y)\f^2 x
\end{array}
\end{equation*}
and the following identities are valid:
\begin{equation}\label{curSl}
\begin{array}{ll}
\n_x \xi=\f x, \qquad &\left(\n_x \eta \right)(y)=g(x,\f y),\\
R(x,y)\xi=-\eta(y)x+\eta(x)y, \qquad &R(\xi,y)\xi=\f^2y, \\
\rho(x,\xi)=-2n\, \eta(x),\qquad 				&\rho(\xi,\xi)=-2n,
\end{array}
\end{equation}
where {$R$ and $\rho$ denote} the curvature tensor and the Ricci tensor, respectively. Let us remark that further $\tau$ and $\tilde\tau$ stand for the scalar curvatures on $\M$ regarding $g$ and $\g$, respectively.

%On an arbitrary Riemannian $\Pi$-Manifolds $\M$, there exists a symmetric $(0,2)$-tensor $\rho^*(y,z)=g^{ij}R(e_i,y,z,\f e_j)$ associated with $\rho$ regarding $\f$.
%In the case of para-Sasaki-like $\M$ the following property is valid \cite{IvMaMa2}:
%\begin{equation*}\label{curf-Sl}
%\begin{array}{l}
%R(x,y,\f z,w)-R(x,y,z,\f w)\\[4pt]
%=-\left\{g(y,z)-2\eta(y)\eta(z)\right\}g(x,\f w)
%-\left\{g(y,w)-2\eta(y)\eta(w)\right\}g(x,\f z)\\[4pt]
%+\left\{g(x,z)-2\eta(x)\eta(z)\right\}g(y,\f
%w)+\left\{g(x,w)-2\eta(x)\eta(w)\right\}g(y,\f z).
%\end{array}
%\end{equation*}
%Then, taking the trace of the latter equality for $x=e_i$ and $w=e_j$, we obtain that the tensors $\rho^*$ and $\rho$ are related as follows
%\begin{equation}\label{rho*-Sl}
%\rho^*(y,z)=\rho(y,\f z)+(2n-1)g(y,\f z).
%\end{equation}

In \cite{HM19}, it is proved that for a $(2n+1)$-dimensional para-Sasaki-like Riemannian $\Pi$-manifold $\M$ the following properties of the Ricci operator $Q$ are valid
\begin{gather}\label{nQxiQ-aSl}
(\n_x Q)\xi=-Q\f x+2n\,\f x,
\\[4pt]
\label{nxiQ=Q-aSl}
(\n_\xi Q)y=-2Q\f y.
%\\[4pt]
%\label{etaQxi-aSl}
%\begin{array}{l}
%\eta\bigl((\n_x Q)\xi\bigr)=0,\qquad
%\eta\bigl((\n_\xi Q)y\bigr)=0.
%\end{array}
\end{gather}
%\begin{proof}
%Taking into account \eqref{curSl}, then $Q\xi=-2n\,\xi$ and $\n_x \xi=\f x$ hold. By virtue of them, we obtain immediately the covariant derivative in \eqref{nQxiQ-aSl}.
%
%After that, we apply $\n_z$ to the expression of $R(x,y)\xi$ in \eqref{curSl}
%and using the form of $\n\eta$ in \eqref{curSl}, we get that
%\[
%\left(\n_z R\right)(x,y)\xi=-R(x,y)\f z + g(y,\f z)x - g(x,\f z)y.
%\]
%Taking the trace of the latter equality for $z=e_i$ and $x=e_j$ and using \eqref{rho*-Sl}, we obtain
%\begin{equation*}\label{tr12}
%g^{ij} (\n_{e_i} R)(e_j,y)\xi=Q\f y + 2n \f y.
%\end{equation*}
%By virtue of the following consequence of the second Bianchi identity
%\[
%g^{ij} (\n_{e_i} R)(\xi,y)e_j=\left(\n_y Q\right)\xi-\left(\n_\xi Q\right)y
%\]
%and using the symmetries of $R$ and \eqref{nQxiQ-aSl}, we prove the truthfulness of \eqref{nxiQ=Q-aSl}.
%
%As consequences of \eqref{nQxiQ-aSl} and \eqref{nxiQ=Q-aSl}, we get \eqref{etaQxi-aSl}.
%\end{proof}

According to \cite{HM17}, a manifold $\M$ is said to be \emph{para-Ein\-stein-like} with constants $(a,b,c)$ if $\rho$ satisfies:
\begin{equation*}\label{defEl}
\begin{array}{l}
\rho=a\,g +b\,\g +c\,\eta\otimes\eta.
\end{array}
\end{equation*}
For $b=0$ or $b=c=0$, the manifold is called an \emph{$\eta$-Einstein manifold} or an \emph{Einstein manifold}, respectively. If $a$, $b$, $c$ are functions on $\mathcal{M}$, then $\M$ is called an \emph{almost para-Einstein-like}, an \emph{almost $\eta$-Einstein manifold} or an \emph{almost Einstein manifold} in the respective cases.

%Let is consider a $(2n+1)$-dimensional Riemannian $\Pi$-manifold $\M$ which is para-Sasaki-like and para-Einstein-like with constants $(a,b,c)$. Tracing \eqref{defEl} and using the last equalities of \eqref{curSl}, we have: \cite{HM17}
%\begin{equation}\label{tauElSl2}
%a+b+c=-2n,\qquad \tau=2n(a-1),
%\end{equation}
%where $\tau$ stands for the scalar curvature with respect to $g$ of $\M$.
%Moreover, for the scalar curvature $\tilde\tau$ with respect to $\g$ on $\M$ we obtain
%\begin{equation}\label{abctau*-ElSl}
%\tilde\tau=2n(b-1).
%\end{equation}
%
%Taking into account \eqref{tauElSl2} and \eqref{abctau*-ElSl}, expression \eqref{defEl} gets the following form
%\begin{equation*}\label{defElSl}
%\begin{array}{l}
%\rho=\left(\dfrac{\tau}{2n}+1\right)g +\left(\dfrac{\tilde\tau}{2n}+1\right)\g
%+\left(-2(n+1)-\dfrac{\tau+\tilde\tau}{2n}\right)\eta\otimes \eta.
%\end{array}
%\end{equation*}

Let us recall from \cite{HM18} the following
\begin{proposition}[\cite{HM18}]\label{prop:El-Dtau}
Let $\M$ be a $(2n+1)$-dimensional para-Sasaki-like Riemannian $\Pi$-manifold. If $\M$ is almost para-Einstein-like  with functions $(a,b,c)$
then the scalar curvatures $\tau$ and $\tilde\tau$ are constants
\begin{equation*}\label{El-Dtauxi}
\tau = const, \qquad \tilde\tau=-2n
\end{equation*}
and $\M$ is $\eta$-Einstein with constants
\[
(a,b,c)=\left(\frac{\tau}{2n}+1,\,0,\,-2n-1-\frac{\tau}{2n}\right).
\]
\end{proposition}

\subsection{Para-Ricci-like soliton on Riemannian $\Pi$-manifold}
A para-Ricci-like soliton with potential $v$ and constants $(\lm,\mu,\nu)$ is determined by the following identity for its Ricci tensor $\rho$:
\begin{equation}\label{defRl}
\begin{array}{l}
\rho=-\frac12 \mathcal{L}_{v} g - \lm\, g - \mu\, \g - \nu\, \eta\otimes \eta,
\end{array}
\end{equation}
where $\mathcal{L}$ stands for the Lie derivative.
In \cite{HM17}, it is considered the soliton potential to be the Reeb vector field, \ie $v=\xi$. In \cite{HM18} and \cite{HM19}, this notion is generalized and there are investigated two cases for the soliton's potential -- when it is pointwise collinear with $\xi$, \ie $v=k\,\xi$ for a differentiable function $k$ on $\mathcal{M}$, and when it is an arbitrary vector field $v$, respectively.

If $\mu=0$ or $\mu=\nu=0$, then \eqref{defRl} defines an $\eta$-Ricci soliton or a Ricci soliton on $\M$, respectively. In the case when $\lm$, $\mu$, $\nu$ are functions on $\mathcal{M}$, the respective soliton is called almost para-Ricci-like soliton, almost $\eta$-Ricci soliton or almost Ricci soliton.

%Let us consider the latter type, \ie a Riemannian $\Pi$-manifold $\M$ admits a para-Ricci-like soliton with arbitrary potential vector field $v$ and constants $(\lm,\mu,\nu)$. Then, its Ricci tensor $\rho$ satisfies the following:
%\begin{equation}\label{defRl-v}
%\begin{array}{l}
%\rho=-\frac12 \mathcal{L}_{v} g - \lm\, g - \mu\, \g - \nu\, \eta\otimes \eta.
%\end{array}
%\end{equation}

Now, we consider a para-Ricci-like soliton with potential which is a gradient of a differentiable function.
\begin{definition}
A Riemannian $\Pi$-manifold $\M$ admits a gradient almost para-Ricci-like soliton with potential vector field $v$ which is a gradient of a differentiable function $f$, i.e. $v=\grad f$, and functions $(\lm,\mu,\nu)$ on $\mathcal{M}$ if its Ricci tensor $\rho$ satisfies the following:
\begin{equation}\label{Hess}
\rho=-\Hess f-\lm g-\mu\g -\nu\eta\otimes\eta,
\end{equation}
where $\Hess$ stands for the Hessian operator with respect to $g$, i.e.
\begin{equation}\label{Hess2}
(\Hess f)(x,y) = (\n_x \D f)(y)=g(\n_x \grad f,y).
\end{equation}
\end{definition}

%}
%
%\subsection{Gradient almost Ricci-like solitons}
%
%{\color{red}

%Let $\lm$, $\mu$, $\nu$ in \eqref{defRl-v} be functions on $\mathcal{M}$.
%If its potential $v$ is a gradient of a differentiable function $f$, i.e. $v=\grad f$, then the soliton is called a \emph{gradient almost Ricci-like soliton} of $\M$. In this case \eqref{defRl-v} is reduced to the following condition
%\begin{equation}\label{Hess}
%\Hess f+\rho+\lm g+\mu\g +\nu\eta\otimes\eta=0,
%\end{equation}
%where $\Hess$ denotes the Hessian operator with respect to $g$, i.e. $\Hess f$ is defined by
%\begin{equation}\label{Hess2}
%(\Hess f)(x,y) := (\n_x \D f)(y)=g(\n_x \grad f,y).
%\end{equation}
%% -- symmetric.
%Taking the trace of \eqref{Hess}, we obtain
% \begin{equation*}\label{Lap}
%\Delta f+\tau+(2n+1)\lm+\mu+\nu=0,
%\end{equation*}
%where $\Delta:=\tr\circ\Hess$ is the Laplacian operator of $g$.
%Also for the Laplacian of $f$, the formula $\Delta f = \Div(\grad f)$ is valid, where $\Div$ stands for the divergence operator.

Further, we omit the trivial case of the gradient para-Ricci-like soliton when $f$ is constant.

\section{Proof of the main results}

Using \eqref{Hess} and \eqref{Hess2}, we get the following
\begin{equation}\label{nv}
\n_x v=-Qx-\lm x-\mu \f x-(\mu+\nu)\eta(x)\xi.
\end{equation}

Taking into acount \eqref{nv}, we calculate that
\begin{equation}\label{Rxyv-aSl}
\begin{array}{l}
R(x,y)v=-\left(\n_x Q\right)y+\left(\n_y Q\right)x\\[4pt]
\phantom{R(x,y)v=}
+\{\D\lm(y)+\mu\eta(y)\}x-\{\D\lm(x)+\mu\eta(x)\}y\\[4pt]
\phantom{R(x,y)v=}
+\{\D\mu(y)+(\mu+\nu)\eta(y)\}\f x-\{\D\mu(x)+(\mu+\nu)\eta(x)\}\f y\\[4pt]
\phantom{R(x,y)v=}
+\D(\mu+\nu)(y)\eta(x) \xi- \D(\mu+\nu)(x)\eta(y) \xi,
\end{array}
\end{equation}
which for $x=\xi$ takes the following form
\begin{equation*}\label{R1-}
\begin{array}{l}
R(\xi,y)v=-\left(\n_\xi Q\right)y+\left(\n_y Q\right)\xi+\{\D\lm(\xi)+\mu\}\f^2 y-\{\D\mu(\xi)+\mu+\nu\}\f y\\[4pt]
\phantom{R(\xi,y)v=}
+\D(\lm+\mu+\nu)(y)\xi-\D(\lm+\mu+\nu)(\xi)\eta(y) \xi.
\end{array}
\end{equation*}
Applying \eqref{nQxiQ-aSl} and \eqref{nxiQ=Q-aSl} in the latter equality, we obtain
\begin{equation}\label{R1}
\begin{array}{l}
R(\xi,y)v=Q\f y +\{\D\lm(\xi)+\mu\}\f^2y-\{\D\mu(\xi)+\mu+\nu-2n\}\f y\\[4pt]
\phantom{R(\xi,y)v=3Q\f y}
+\D(\lm+\mu+\nu)(y)\xi-\D(\lm+\mu+\nu)(\xi)\eta(y) \xi.
\end{array}
\end{equation}

On the other hand, if we substitute $z=v$ in the equality for $R(\xi,y)z$ in \eqref{curSl}, we get the following expression %in this case
\begin{equation}\label{Rxiyv-aSl}
R(\xi,y)v=-\D f(y)\xi+\D f(\xi)y.
\end{equation}

Therefore, \eqref{R1} and \eqref{Rxiyv-aSl} deduce that
\begin{equation}\label{nxiQ-aSl}
\begin{array}{l}
Q\f y=-\{\D(\lm-f)(\xi)+\mu\}\f^2y+\{\D\mu(\xi)+\mu+\nu-2n\}\f y\\[4pt]
\phantom{Q\f y=}
-\D(\lm+\mu+\nu+f)(y)\xi+\D(\lm+\mu+\nu+f)(\xi)\eta(y) \xi.
\end{array}
\end{equation}

Taking into account the property $Q\circ\f=\f\circ Q$ for a para-Sasaki-like manifold, equality \eqref{nxiQ-aSl} implies the following
\begin{equation}\label{lmnf+=-aSl}
\D(\lm+\mu+\nu+f)(y)=\D(\lm+\mu+\nu+f)(\xi)\eta(y).
\end{equation}
Therefore, \eqref{nxiQ-aSl} and \eqref{R1} take the following form
\begin{equation*}\label{Qfi-aSl}
Q\f y=-\{\D(\lm-f)(\xi)+\mu\}\f^2y+\{\D\mu(\xi)+\mu+\nu-2n\}\f y,
\end{equation*}
\begin{equation*}\label{R1=}
\begin{array}{l}
R(\xi,y)v=Q\f y +\{\D\lm(\xi)+\mu\}\f^2y-\{\D\mu(\xi)+\mu+\nu-2n\}\f y\\[4pt]
\phantom{R(\xi,y)v=Q\f y}
-\{\D f(y)-\D f(\xi)\eta(y)\}\xi.
\end{array}
\end{equation*}
Then, we have
\begin{equation}\label{Rxiyvz}
\begin{array}{l}
R(\xi,y,v,z)\,{=} \rho(y,\f z) +\{\D\lm(\xi)+\mu\}g(\f y,\f z)
\\[4pt]
\phantom{R(\xi,y,v,z)\,{=}}
- \{\D\mu(\xi)+\mu+\nu-2n\}g( y,\f z)\\[4pt]
\phantom{R(\xi,y,v,z)\,{=}}
-\{\D f(y)-\D f(\xi)\eta(y)\}\eta(z).
\end{array}
\end{equation}
%
%
%
%Combining \eqref{R1} and \eqref{Rxiyv-aSl}, we find the following formula
%\begin{equation}\label{R1==}
%\begin{array}{l}
%3Q\f y= -\mu \f^2 y +\D(\lm-f)(\xi)y+\{\D\mu(\xi)+\mu+\nu+2n\}\f y\\[4pt]
%\phantom{3Q\f y=}
%-\D(\lm+\mu+\nu-f)(y)\xi+\D(\mu+\nu)(\xi)\eta(y) \xi.
%\end{array}
%\end{equation}

From another point of view, the form of $R(x,y)\xi$ from \eqref{curSl} and equality \eqref{Rxyv-aSl} imply the following
\begin{gather}
R(x,y,\xi,v)=-\D f(x)\eta(y)+\D f(y)\eta(x),\label{R3}\nonumber
\\[4pt]
\begin{array}{l}\label{R4}
R(x,y,v,\xi)=-\eta\left((\n_x Q)y-(\n_y Q)x\right)+\D(\lm+\mu+\nu)(y)\eta(x)\nonumber\\[4pt]
\phantom{R(x,y,v,\xi)=-\eta\left((\n_x Q)y-(\n_y Q)x\right)}
-\D(\lm+\mu+\nu)(x)\eta(y).
\end{array}
\end{gather}
Combining the latter two equalities, we get that
\begin{equation*}\label{etaQ-aSl}
\begin{array}{l}
\eta\bigl((\n_x Q)y-(\n_y Q)x\bigr)=-\,\D(\lm+\mu+\nu+f)(x)\eta(y)\\[4pt]
\phantom{\eta\left((\n_x Q)y-(\n_y Q)x\right)=}
+\D(\lm+\mu+\nu+f)(y)\eta(x),
\end{array}
\end{equation*}
which, taking into account \eqref{lmnf+=-aSl}, is simplified to
\begin{equation*}\label{etaQ=00-aSl}
\eta\bigl((\n_x Q)y-(\n_y Q)x\bigr)=0.
\end{equation*}

%Then, equalities \eqref{etaQ=00-aSl}, \eqref{etaQ-aSl}, \eqref{R4} and \eqref{R3} imply
%\[
%\D f(x)\eta(y)-\D f(y)\eta(x)=0
%\]
%and its consequence
%\begin{equation}\label{Df}
%\D f(x)=\D f(\xi)\eta(x).
%\end{equation}
%
%Substituting the latter equality in \eqref {R1=}, gives us
%\[
%R(\xi,y)v=-Q\f y +\{\D\lm(\xi)+\mu\}\f^2y-\{\D\mu(\xi)+\mu+\nu+2n\}\f y
%\]
%and therefore we have
%\begin{equation}\label{Rxiyvz}
%\begin{array}{l}
%R(\xi,y,v,z)\,{=} - \rho(y,\f z) +\{\D\lm(\xi)+\mu\}g(\f y,\f z)
%\\[4pt]
%\phantom{R(\xi,y,v,z)\,{=}}- \{\D\mu(\xi)+\mu+\nu+2n\}g( y,\f z).
%\end{array}
%\end{equation}

On the other hand, the expression of $R(\xi,y)z$ from \eqref{curSl} % and \eqref{Df}
gives
\[
R(\xi,y,z,v)\,{=}-\D f(\xi)g(\f y, \f z)+\{\D f(y)-\D f(\xi)\eta(y) \}\eta(z).
\]
The latter result, \eqref{Rxiyvz} and $\rho(x,\xi)$ from \eqref{curSl} imply
\begin{equation*}\label{rho-v}
\begin{array}{l}
\rho(y,z)=\{\D\mu(\xi)+\mu+\nu-2n\}g(\f y, \f z)-\{\D(\lm-f)(\xi)+\mu\}g(y,\f z)\\[4pt]
\phantom{\rho(y,z)=}
-2n\eta(y)\eta(z),
\end{array}
\end{equation*}
which can be rewritten in the following form
\[
%\begin{equation}\label{rho-v}
\begin{array}{l}
\rho=\,\{\D\mu(\xi)+\mu+\nu-2n\}g-\{\D(\lm-f)(\xi)+\mu\}\g\\[4pt]
\phantom{\rho=}
-\{\nu-\D(\lm-\mu-f)(\xi)\}\eta\otimes\eta.
\end{array}
\]
Therefore, $\M$ is almost para-Einstein-like with coefficient functions
\begin{equation}\label{abc-aElSlRl}
\begin{array}{c}
a\,{=}\D\mu(\xi)+\mu+\nu-2n,\qquad
b=-\D(\lm-f)(\xi)-\mu,\\[4pt]
c\,{=}\D(\lm-\mu-f)(\xi)-\nu.
\end{array}
\end{equation}
Then, using \propref{prop:El-Dtau} and \eqref{abc-aElSlRl}, we have that
\begin{equation}\label{prop-1}
a\,{=}\frac{\tau}{2n}+1,\qquad b\,{=}0,\qquad c\,{=}-2n-1-\frac{\tau}{2n},
\end{equation}
\begin{equation}\label{prop-2}
\tau\,{=}-2n\{\D\mu(\xi)+\mu+\nu-2n-1\},\qquad \ttt\,{=}-2n.
\end{equation}

Taking into account \eqref{prop-2}, we get $\D\ttt=0$ and using $\D\ttt(y)=2(\tau+2n)\eta(y)$, we obtain
\[
\tau=-2n.
\]
Therefore, substituting the latter result in \eqref{prop-1}, we get
\[
(a,b,c)=(0,0,-2n),
\]
which completes the proof.

\section{Example}\label{Ex2}

Let us consider the following example given in \cite{HM19}.
Let $\mathcal{M}$ be a set of points in $\R^3$ with coordinates $(x^1,x^2,x^3)$, which is equipped with a Riemannian $\Pi$-structure $(\f, \xi, \eta, g)$ as follows
\begin{equation*}\label{strEx1-loc}
\begin{array}{c}
g\left(\DD_1,\DD_1\right)=g\left(\DD_2,\DD_2\right)=\cosh 2x^3,\qquad
g\left(\DD_1,\DD_2\right)=\sinh 2x^3,
\\[4pt]
g(\DD_1,\DD_3)=g(\DD_2,\DD_3)=0,\qquad g(\DD_3,\DD_3)=1,
\\[4pt]
\f  \DD_1=\DD_2,\qquad \f \DD_2=\DD_1,\qquad \xi=\DD_3,
\end{array}
\end{equation*}
where by $\DD_1$, $\DD_2$, $\DD_3$ we denote
$\frac{\DD}{\DD{x^1}}$,
$\frac{\DD}{\DD{x^2}}$, $\frac{\DD}{\DD{x^3}}$, respectively.
Then, the triad
\begin{equation}\label{eiloci}
\{e_1,e_2,e_3\}=\{\cosh x^3 \DD_1-\sinh x^3 \DD_2,\,-\sinh x^3 \DD_1+\cosh x^3 \DD_2,\,\DD_3\}
\end{equation}
forms an orthonormal $\f$-basis of $T_p\mathcal{M}$, $p\in \mathcal{M}$. Therefore, we have
\begin{equation*}\label{strEx1}
\begin{array}{c}
g(e_i,e_i)=1,\quad
g(e_i,e_j)=0,\quad
i,j\in\{1,2,3\},\; i\neq j,
\\[4pt]
\f  e_1=e_2,\qquad \f e_2=e_1,\qquad \xi=e_3.
\end{array}
\end{equation*}

Then, the constructed manifold $\M$ is a 3-dimensional para-Sasaki-like Riemannian $\Pi$-manifold \cite{HM19}. Moreover, in the cited paper for $\M$ is obtained that
\begin{equation}\label{rhotauk}
\tau=\ttt=-2,\qquad \rho=-2\eta\otimes\eta.
\end{equation}

Now, let $f$ be a differentiable function on $\mathcal{M}$, determined by
\begin{equation*}\label{f}
f=\frac{1}{2}p\,\{(x^1)^2+(x^2)^2\} -x^2 + q\,x^3,
\end{equation*}
where $p$ and $q$ are arbitrary constants.
Then, the gradient of $f$ with respect to $g$ has the following form
\begin{equation}\label{gradf}
\begin{array}{l}
\grad f = \,\{p\,x^1\cosh x^3 + (p\,x^2-1)\sinh x^3\}e_1\\[4pt]
\phantom{\grad f =}
+\{p\,x^1\sinh x^3 + (p\,x^2-1)\cosh x^3\}e_2 + q\, e_3.
\end{array}
\end{equation}

Taking into account \eqref{eiloci}, we calculate the  components of
$\LL_{\grad f}g$ as follows
\[
\left(\LL_{\grad f}g\right)_{11}=\left(\LL_{\grad f}g\right)_{22}=-2p,\qquad
\left(\LL_{\grad f}g\right)_{12}=-2q.
\]
Therefore, we have the following expression
\begin{equation}\label{ex-1}
\left(\LL_{\grad f}g\right)=-2p\, g -2q \g +2(p+q)\,\eta\otimes\eta.
\end{equation}

Equality \eqref{ex-1} coincides with
\[
\LL_v g = -2\,c_1\,g   +2(c_2+c_3)\g   +2(c_1-c_2-c_3)\eta\otimes\eta
\]
obtained for the same example in \cite{HM19} when $p=c_1$, $q=-(c_2+c_3)$.
Therefore, the manifold $\M$ admits a gradient almost para-Ricci-like soliton with potential $v=\grad f$ determined by \eqref{gradf}
and $(\lm,\mu,\nu)$ are the following
\[
\lm=p,\qquad
\mu=q,\qquad
\nu=-p-q+2.
\]

In conclusion, the constructed 3-dimensional example of a para-Sasaki-like Riemannian $\Pi$-manifold $\M$
with \eqref{rhotauk} and gradient para-Ricci-like soliton supports \thmref{thm:grad}.

%%%%%%%%%%%%%%%%%%%%%%%%%%%%%%%%%%%%%%%%%%

\end{document}